\documentclass[a4paper,english,oneside]{article}
\usepackage[T1]{fontenc}
\usepackage[utf8]{inputenc} 

\usepackage{geometry}
\usepackage{amsmath}
\usepackage{amsfonts}
\usepackage{amssymb}
\usepackage{amsthm}
\usepackage[english]{babel}
\usepackage{color}
\usepackage{lmodern}
\usepackage{pstricks-add}
\usepackage{hyperref}
\hypersetup{hidelinks}

\usepackage{mathscinet}
\usepackage{enumitem}
\usepackage[cal=boondoxo,bb=ams]{mathalfa}
\usepackage{microtype}

\usepackage{mathtools}
\usepackage{esint}

\theoremstyle{plain}
\newtheorem{theorem}{Theorem}[section]

\newtheorem{proposition}[theorem]{Proposition}

\theoremstyle{definition}

\theoremstyle{remark}
\newtheorem{remark}{Remark}

    \DeclareMathOperator\supp{supp}
    \DeclareMathOperator\meas{meas}

    \DeclareMathOperator\Fuj{Fuj}
    \DeclareMathOperator\Str{Str}
    \DeclareMathOperator\Gla{Gla}

\begin{document}

\title{A blow-up result for a semilinear wave equation with scale-invariant damping and mass and nonlinearity of derivative type}

\author{Alessandro Palmieri$^\mathrm{a}$, Ziheng Tu$^\mathrm{b}$}

\date{\small{$^\mathrm{a}$ Department of Mathematics, University of Pisa,  Largo B. Pontecorvo 5, 56100 Pisa, Italy}\\  \small{$^\mathrm{b}$ School of Data Science, Zhejiang University of Finance and Economics, 310018 Hangzhou, China}\\[2ex] \normalsize{\today} }

\maketitle

\begin{abstract}

In this note, we prove blow-up results for semilinear wave models with damping and mass in the scale-invariant case and with nonlinear terms of derivative type. We consider the single equation and the weakly coupled system. In the first case we get a blow-up result for exponents below a certain shift of the Glassey exponent. For the weakly coupled system we find as critical curve a shift of the corresponding curve for the weakly coupled system of semilinear wave equations with the same kind of nonlinearities. Our approach follows the one for the respective classical wave equation by Zhou Yi. In particular, an explicit integral representation formula for a solution of the corresponding linear scale-invariant wave equation, which is derived by using Yagdjian's integral transform approach, is employed in the blow-up argument. While in the case of the single equation we may use a comparison argument, for the weakly coupled system an iteration argument is applied.

\end{abstract}

\begin{flushleft}
\end{flushleft}

\begin{flushleft}
\textbf{AMS Classification (2010)} Primary:  35B44, 35L71; Secondary: 35B33, 35C15
\end{flushleft}

\section{Introduction}

In this work we prove a blow-up result for the semilinear wave equation with time-dependent damping and mass in the scale-invariant case and with nonlinearity of derivative type, namely,
\begin{align}\label{semilinear CP derivative type}
\begin{cases} u_{tt}-\Delta u +\frac{\mu}{1+t}u_t+\frac{\nu^2}{(1+t)^2}u=|\partial_t u|^p, &  x\in \mathbb{R}^n, \ t>0,\\
u(0,x)=\varepsilon u_0(x), & x\in \mathbb{R}^n, \\ u_t(0,x)=\varepsilon u_1(x), & x\in \mathbb{R}^n,
\end{cases}
\end{align} where $\mu,\nu^2$ are nonnegative constants, $p>1$ and $\varepsilon$ is a positive constant describing the smallness of Cauchy data. Let us introduce the quantity
\begin{align}
\delta\doteq (\mu-1)^2-4 \nu^2. \label{def delta}
\end{align} Recently, semilinear wave equation with scale-invariant damping and mass terms and power nonlinearity $|u|^p$ has been studied in several papers. It turns out that if $\delta$ is ``large'', that is, for $\delta\geq (n+1)^2$, the critical exponent for 
\begin{align}\label{semilinear CP power nonlinear}
\begin{cases} u_{tt}-\Delta u +\frac{\mu}{1+t}u_t+\frac{\nu^2}{(1+t)^2}u=| u|^p, &  x\in \mathbb{R}^n, \ t>0,\\
u(0,x)=\varepsilon u_0(x), & x\in \mathbb{R}^n, \\ u_t(0,x)=\varepsilon u_1(x), & x\in \mathbb{R}^n,
\end{cases}
\end{align} is given by the shift $p_{\Fuj}\left(n+\frac{\mu-1}{2}-\frac{\sqrt{\delta}}{2}\right)$ of the \emph{Fujita exponent} $p_{\Fuj}(n)\doteq 1+\frac{2}{n}$ (cf. \cite{Wak14,Abb15,NPR16,PalRei17,Pal17}). This follows from the fact that the critical exponent of the semilinear classical damped wave equation with power nonlinearity is the Fujita exponent and for for large $\delta$ \eqref{semilinear CP power nonlinear} has simliar properties to this model somehow. On the other hand, for ``small'' and nonnegative value of $\delta$ the critical exponent for \eqref{semilinear CP power nonlinear} should reasonably be the shift $p_{\Str}(n+\mu)$ of the \emph{Strauss exponent} $p_{\Str}(n)$ the critical exponent for the semilinear wave equations with power nonlinearity (named after the author of \cite{Str81}, where a conjecture for the critical exponent for the semilinear wave equation with $|u|^p$ as nonlinear term is done), which is the positive root of the quadratic equation $(n-1)p^2-(n+1)p-2=0$ (cf. \cite{John79,Kato80,Glas81B,Sid84,Scha85pn,Lin90,Zhou92,Zhou93,YZ06,Zhou07,TakWak11,ZH14} for the necessity part and \cite{John79,Glas81,LinSog95,LinSog96,Geo97,Tat01,Jiao03} for the sufficiency part or \cite{Ku94,KuKu96} for the radial symmetric case). This conjecture for the scale-invariant model is still open for the sufficiency part (for the necessity part, that is the blow-up results, see \cite{DLR15,Wak16,LTW17,IS17,TL1709,PalRei18,PT18,KatoSak18}), even though some partial results in the special case $\delta=1$ have been proved for $n\geq 3$ in the radial symmetric case (see \cite{DabbLuc15,Pal18odd} for the odd dimensional case and \cite{Pal18even} for the even dimensional case, respectively). This peculiarity of a ``parabolic-like'' behavior for large values of $\delta$ and of ``wave-like'' behavior for small values of $\delta$ has been showed also for the corresponding weakly coupled system (cf. \cite{ChenPal18,Pal19}).

In the case of the Cauchy problem for the semilinear wave equation with nonlinearity of derivative type 
\begin{align}\label{semilinear CP classical wave derivative type}
\begin{cases} u_{tt}-\Delta u =|\partial_t u|^p, &  x\in \mathbb{R}^n, \ t>0,\\
u(0,x)=\varepsilon u_0(x), & x\in \mathbb{R}^n, \\ u_t(0,x)=\varepsilon u_1(x), & x\in \mathbb{R}^n,
\end{cases}
\end{align} the critical exponent is the so-called \emph{Glassey exponent} $p_{\Gla}(n)\doteq \frac{n+1}{n-1}$.
We refer to the classical works \cite{Joh81,Sid83,Mas83,Sch85,Ram87,Age91,HT95,Tzv98,Zhou01,HWY12} for the proof of this conjecture, although up to the knowledge of the author the global existence in the supercritical case for the not radial symmetric case in high dimensions is still open. Recently, in \cite{LT18Glass} a blow-up result for $1<p\leq p_{\Gla}(n)$ has been proved for a semilinear damped wave model in the scattering case, that is, when the time-dependent coefficient of the damping term $b(t)u_t$ is nonnegative and summable.

Therefore, according to what happens for the semilinear Cauchy problem \eqref{semilinear CP power nonlinear} it would be natural to find as critical exponent for \eqref{semilinear CP derivative type} a suitable shift for the Glassey exponent. Purpose of this paper is to prove under certain sign assumptions for the Cauchy data a blow-up result for the Cauchy problem \eqref{semilinear CP derivative type} provided that the exponent in the nonlinear term satisfies $1<p\leq p_{\Gla}(n+\sigma)$ with the shift $\sigma$ defined by  
\begin{align}\label{def shift sigma}
\sigma\doteq \begin{cases} \mu +1-\sqrt{\delta} & \mbox{if } \ \delta \in [0,1), \\ \mu  & \mbox{if } \ \delta\geq 1. \end{cases}
\end{align}
 As byproduct of the comparison argument that will be employed we get an upper bound estimate for the lifespan in terms of $\varepsilon$ as well.

Let us consider the weakly coupled system of semilinear wave equations with damping and mass in the scale-invariant case and nonlinearities of derivative type, that is,
\begin{align}\label{semilinear WCS derivative type}
\begin{cases} u_{tt}-\Delta u +\frac{\mu_1}{1+t}u_t+\frac{\nu^2_1}{(1+t)^2}u=|\partial_t v|^p, &  x\in \mathbb{R}^n, \ t>0,\\
 v_{tt}-\Delta v +\frac{\mu_2}{1+t}v_t+\frac{\nu^2_2}{(1+t)^2}v=|\partial_t u|^q, &  x\in \mathbb{R}^n, \ t>0,\\
u(0,x)=\varepsilon u_0(x), \ \ v(0,x)=\varepsilon v_0(x), & x\in \mathbb{R}^n, \\ u_t(0,x)=\varepsilon u_1(x), \ \  v_t(0,x)=\varepsilon v_1(x), & x\in \mathbb{R}^n,
\end{cases}
\end{align} where $\mu_1,\mu_2,\nu^2_1,\nu^2_2$ are nonnegative constants, $p,q>1$ and $\varepsilon$ is a positive constant describing the smallness of Cauchy data. Similarly to the case of a single equation, we introduce the quantities $\delta_j\doteq (\mu_j-1)^2-4\nu^2_j$ for $j=1,2$.  

The machinery, that we are going to develop in the case of a single semilinear equation, works nicely also in the case of this weakly coupled system.

In order to understand our blow-up result for \eqref{semilinear WCS derivative type}, we shall first recall some results which are known in the literature for the semlinear weakly coupled system of wave equation with nonlinear terms of derivative type, namely,
\begin{align}\label{WCS classical}
\begin{cases} u_{tt}-\Delta u =|\partial_t v|^p, &  x\in \mathbb{R}^n, \ t>0,\\
 v_{tt}-\Delta v =|\partial_t u|^q, &  x\in \mathbb{R}^n, \ t>0,\\
u(0,x)=\varepsilon u_0(x), \ \ v(0,x)=\varepsilon v_0(x), & x\in \mathbb{R}^n, \\ u_t(0,x)=\varepsilon u_1(x), \ \  v_t(0,x)=\varepsilon v_1(x), & x\in \mathbb{R}^n,
\end{cases}
\end{align}

The non-exitence of global in time solutions to \eqref{WCS classical} (which corresponds to \eqref{semilinear WCS derivative type} in the case $\mu_1= \mu_2=0$ and $\nu_1^2=\nu_2^2=0$) has been studied in \cite{Deng99,Xu04}, while the existence part has been proved in the three dimensional and radial case in \cite{KKS06}.  Recently, in \cite[Section 8]{ISW18} the upper bound for the lifespan has been derived. 
Summarizing the main results of these works we can see that \begin{align}\label{critical case classical model WCS}
\max\{\Lambda(n,p,q), \Lambda(n,q,p)\}=0,
\end{align} is the critical line in the $(p,q)$-plane for the semilinear weakly coupled system \eqref{WCS classical} where
\begin{align}\label{def Lambda}
\Lambda(n,p,q)\doteq \frac{p+1}{pq-1}-\frac{n-1}{2}.
\end{align} 

 Let us recall the meaning of critical curve for a weakly coupled system: if the exponents $p,q>1$ satisfy $\max\{\Lambda(n,p,q), \Lambda(n,q,p)\}<0$ (\emph{supercritical case}), then, it is possible to prove a global existence result for small data solutions; on the contrary, for $\max\{\Lambda(n,p,q), \Lambda(n,q,p)\}\geq 0$ it is possible to prove the nonexistence of global in time solutions regardless the smallness of the Cauchy data and under certain sign assumptions for them. Let us point out that, according to the results we quoted above, the conjecture that the critical line for \eqref{WCS classical} is given by \eqref{critical case classical model WCS} has be shown to be true only partially, as the global existence of small data solutions has been proved only in the 3-dimensional and radial symmetric case. 

In the massless case ($\nu_1^2=\nu_2^2=0$) and \emph{scattering producing} case, that is, if we consider time-dependent, nonnegative and summable coefficients $b_1(t),b_2(t)$ instead of $\mu_1(1+t)^{-1},\mu_2(1+t)^{-1}$, really recently in \cite{PalTak19dt} a blow-up result has been proved in the same range for the exponents $(p,q)$ as for the corresponding not damped case, namely, for $(p,q)$ such that $\max\{\Lambda(n,p,q), \Lambda(n,q,p)\}\geq 0$ is satisfied.

Consequently, coming back to the weakly coupled system in the scale-invariant case \eqref{semilinear WCS derivative type}, we may expect as critical curve in the $(p,q)$-plane a curve with branches that are shifts of the branches of the critical curve for \eqref{WCS classical}. Indeed, due to the blow-up result for \eqref{semilinear WCS derivative type} which we are going to state in the next section, we may conjecture $$\max\{\Lambda(n+\sigma_1,p,q), \Lambda(n+\sigma_2,q,p)\}=0$$ as critical curve, where $\sigma_1,\sigma_2$ are defined analogously as in \eqref{def shift sigma}.

\section{Main results}

\begin{theorem} \label{Thm blow-up derivative type nonlinearity} Let $n\geq 1$ and let $\mu,\nu^2$ be nonnegative constants such that $\delta\geq 0$. We consider $u_0,u_1$ compactly supported in $B_R$ such that $u_0$ and $u_1$ are nonnegative functions if $\delta\geq 1$; else, if $\delta\in [0,1)$ we assume that $u_0=0$ and $u_1$ is a nonnegative function.

 Let $1<p\leq p_{\Gla}(n+\sigma)$ be the exponent of the nonlinearity of derive type, where $\sigma$ is defined by \eqref{def shift sigma}. Then, there exists $\varepsilon_0=\varepsilon_0(n,p,\mu,\nu^2,u_0,u_1,R)>0$ such that for any $\varepsilon\in (0,\varepsilon_0]$ if $u$ is a local in time solution to \eqref{semilinear CP derivative type}, 
$u$ blows up in finite time. Furthermore, the following upper bound estimate for the lifespan of the solution holds
\begin{align} \label{upper bound lifespan semilinear eq}
T(\varepsilon) \leq \begin{cases} C \varepsilon^{-\left(\frac{1}{p-1}-\frac{n+\sigma-1}{2}\right)^{-1}} & \mbox{if} \ \ 1<p<p_{\Gla}(n+\sigma), \\ \exp\big(C \varepsilon^{-(p-1)}\big) & \mbox{if} \ \ p=p_{\Gla}(n+\sigma), \end{cases}
\end{align} where the positive constant $C$ is independent of $\varepsilon$.
\end{theorem}

\begin{remark} The sign assumptions on Cauchy data in the statement of Theorem \ref{Thm blow-up derivative type nonlinearity} are done with the purpose to ensure a suitable control from below for a function which depends on the solution of the homogeneous linear problem related to \eqref{semilinear CP derivative type}. 
\end{remark}

\begin{theorem} \label{Thm blow-up WCS der type nonlinearities} Let $n\geq 1$ and let $\mu_1,\mu_2,\nu^2_1,\nu^2_2$ be nonnegative constants such that $\delta_1,\delta_2\geq 0$. We consider data $u_0,u_1,v_0,v_1$ that are nonnegative and  compactly supported in $B_R$ functions. Moreover, we assume $u_0=0$ $($respectively $v_0=0)$ in the case $\delta_1\in [0,1)$ $($respectively $\delta_2\in~ [0,1))$. Let us assume that the exponents of the nonlinearities of derive type $p,q>1$ satisfy
\begin{align}\label{critical line p,q}
\max\left\{\frac{p+1}{pq-1}-\frac{n+\sigma_1-1}{2},\frac{q+1}{pq-1}-\frac{n+\sigma_2-1}{2}\right\}\geq 0,
\end{align} where 
\begin{align}\label{def shift sigma k}
\sigma_k\doteq \begin{cases} \mu_k +1-\sqrt{\delta_k} & \mbox{if } \ \delta_k \in [0,1), \\ \mu_k  & \mbox{if } \ \delta_k\geq 1, \end{cases} \qquad k=1,2.
\end{align}
 Then, there exists $\varepsilon_0=\varepsilon_0(n,p,q,\mu_1,\mu_2,\nu_1^2,\nu_2^2,u_0,u_1,v_0,v_1,R)>0$ such that for any $\varepsilon\in~ (0,\varepsilon_0]$ if $(u,v)$ is a local in time solution to \eqref{semilinear WCS derivative type}, 
$(u,v)$ blows up in finite time. Furthermore, the following upper bound estimate for the lifespan of the solution holds
\begin{align}\label{upper bound lifespan WCS}
T(\varepsilon) \leq \begin{cases} C \varepsilon^{-\left(\Omega(n,\sigma_1,\sigma_2,p,q)\right)^{-1}} & \mbox{if} \ \ \Omega(n,\sigma_1,\sigma_2,p,q)>0, \\ \exp\Big(C \varepsilon^{-(pq-1)}\Big) \vphantom{\bigg(} & \mbox{if} \ \  \Omega(n,\sigma_1,\sigma_2,p,q)=0 , \\ 
\exp\Big(C \varepsilon^{-\min\left\{\frac{pq-1}{p+1},\frac{pq-1}{q+1}\right\}}\Big) & \mbox{if} \ \  \Lambda(n+\sigma_1,p,q)=\Lambda(n+\sigma_2,q,p)=0, \end{cases}
\end{align} where the positive constant $C$ is independent of $\varepsilon$, $\Lambda(n,p,q)$ is defined by \eqref{def Lambda} and $$\Omega(n,\sigma_1,\sigma_2,p,q)\doteq \max\{\Lambda(n+\sigma_1,p,q),\Lambda(n+\sigma_2,q,p)\}.$$
\end{theorem}

\begin{remark}\label{Remark double critic case}
In the cusp point of the critical line that we found in Theorem \ref{Thm blow-up WCS der type nonlinearities}, that is, for $(p,q)$ such that $\Lambda(n+\sigma_1,p,q)=\Lambda(n+\sigma_2,q,p)=0$, we may specify more explicitly the condition on the lifespan. Indeed, if we denote $$\eta\doteq \frac{n+\sigma_2+1}{2}-\frac{n+\sigma_1-1}{2}q, \quad \xi\doteq \frac{n+\sigma_1+1}{2}-\frac{n+\sigma_2-1}{2}p,$$ then, straightforward computations lead to 
\begin{align*}
p\eta +\xi = (pq-1) \Lambda (n+\sigma_1,p,q) \ \ \mbox{and} \ \ \eta +q\xi = (pq-1) \Lambda (n+\sigma_2,q,p) .
\end{align*} Therefore, on the cusp point of the critical line we get $\eta=\xi=0$ due to the fact that $pq>1$. From $\eta=0$ and $\xi=0$ we obtain the explicit expressions of $(p,q)$ when $\Lambda(n+\sigma_1,p,q)=\Lambda(n+\sigma_2,q,p)=0$, namely,
\begin{align*}
p=\widetilde{p}(n,\sigma_1,\sigma_2)\doteq \frac{n+\sigma_1+1}{n+\sigma_2-1} \  \ \ \mbox{and} \ \ \ q=\widetilde{q}(n,\sigma_1,\sigma_2)\doteq \frac{n+\sigma_2+1}{n+\sigma_1-1}.
\end{align*} Since $\widetilde{p}(n,\sigma_1,\sigma_2)\geq \widetilde{q}(n,\sigma_1,\sigma_2) $ if and only if $\sigma_1\geq \sigma_2$, we can rewrite the last upper bound estimate for the lifespan in \eqref{upper bound lifespan WCS} as follows:
\begin{align*}
T(\varepsilon) \leq \begin{cases}
\exp\Big(C \varepsilon^{-\frac{pq-1}{p+1}}\Big) &  \mbox{if} \ \ \Lambda(n+\sigma_1,p,q)=\Lambda(n+\sigma_2,q,p)=0 \ \ \mbox{and} \ \  \sigma_1\geq \sigma_2, \\
\exp\Big(C \varepsilon^{-\frac{pq-1}{q+1}}\Big) & \mbox{if} \ \ \Lambda(n+\sigma_1,p,q)=\Lambda(n+\sigma_2,q,p)=0 \ \ \mbox{and} \ \  \sigma_1\leq \sigma_2.\end{cases}
\end{align*} Of course, when $\sigma_1=\sigma_2$ these estimates coincide with the estimate for the critical case in \eqref{upper bound lifespan semilinear eq} and $\widetilde{p}(n,\sigma_1,\sigma_1)=\widetilde{q}(n,\sigma_1,\sigma_1)=p_{\Gla}(n+\sigma_1)$.
\end{remark}

Let us illustrate our strategy in the proof of Theorems \ref{Thm blow-up derivative type nonlinearity} and \ref{Thm blow-up WCS der type nonlinearities}: our approach in the proof of the blow-up results is based on the work \cite{Zhou01} for the classical wave equation with nonlinearity of derivative type; therefore, as main tool we need to employ an integral representation formula for the linear and one-dimensional problem associated to \eqref{semilinear CP derivative type}, which generalize d'Alembert's formula in the case of the free wave equation. This formula has been proved really recently in \cite{Pal19RF}. Applying such formula, we end up with a nonlinear ordinary integral inequality (OII) for the single equation \eqref{semilinear CP derivative type} and a system of OIIs for the weakly coupled system \eqref{semilinear WCS derivative type}, respectively. Then, for \eqref{semilinear CP derivative type} a simple comparison argument suffices to prove Theorem \ref{Thm blow-up derivative type nonlinearity}, while in order to prove Theorem \ref{Thm blow-up WCS der type nonlinearities} we shall employ an iteration argument. Furthermore, in the critical case we will combine it with the slicing method.

\section{Proof of Theorem \ref{Thm blow-up derivative type nonlinearity}}

This section is devoted to the proof of Theorem \ref{Thm blow-up derivative type nonlinearity}. Before introducing the suitable function that will allow us to prove the blow-dynamic in the case $1<p\leq p_{\Gla}(n+\sigma)$, we recall the previously mentioned generalization of D'Alembert's  representation formula.

\subsection{Integral representation formula for the 1-dimensional linear case}

In this subsection, we recall a representation formula for the solution of the linear Cauchy problem for a scale-invariant wave equation, namely,
\begin{align}\label{inhomog lin CP n=1}
\begin{cases} u_{tt}- u_{xx} +\frac{\mu}{1+t}u_t+\frac{\nu^2}{(1+t)^2}u=f(t,x), &  x\in \mathbb{R}, \ t>0,\\
u(0,x)=u_0(x), & x\in \mathbb{R}, \\ u_t(0,x)=u_1(x), & x\in \mathbb{R},
\end{cases}
\end{align} where $\mu,\nu^2$ are nonnegative constants. For the proof of this formula one can see \cite[Theorem 1.1]{Pal19RF}.

\begin{proposition} \label{Prop representation formula 1d case}
Let $n=1$ and let $\mu,\nu^2$ be nonnegative constants. Let us assume $u_0\in \mathcal{C}^2(\mathbb{R})$, $u_1\in \mathcal{C}^1(\mathbb{R})$ and $f\in \mathcal{C}^{0,1}_{t,x}([0,\infty)\times \mathbb{R})$. Then, a representation formula for the solution of \eqref{inhomog lin CP n=1} is given by
\begin{align}
u(t,x) &= \frac{1}{2}(1+t)^{-\frac{\mu}{2}}\big(u_0(x+t)+u_0(x-t)\big) \notag\\ & \quad +\frac{1}{2^{\sqrt{\delta}}}\int_{x-t}^{x+t} \Big[ u_0(y) K_0(t,x;y;\mu,\nu^2) + \big(u_1(y)+\mu \, u_0(y)\big) K_1(t,x;y;\mu,\nu^2)\Big]\, \mathrm{d}y \notag \\ & \quad +\frac{1}{2^{\sqrt{\delta}}} \int_0^t \int_{x-t+b}^{x+t-b} f(b,y) E(t,x;b,y;\mu,\nu^2) \, \mathrm{d}y\, \mathrm{d}b, \label{representation formula 1d case}
\end{align} where the kernel functions are defined as follows
\begin{align}
E(t,x;b,y;\mu,\nu^2) & \doteq (1+t)^{-\frac{\mu}{2}+\gamma} (1+b)^{\frac{\mu}{2}+\gamma} \left((t+b+2)^2-(y\!-\!x)^2\right)^{-\gamma} \notag \\  &  \qquad  \times \mathsf{F}\left(\gamma,\gamma;1; \tfrac{(t-b)^2-(y-x)^2}{(t+b+2)^2-(y-x)^2} \right), \label{def E(t,x;b,y)} \\
K_0(t,x;y;\mu,\nu^2) & \doteq -\partial_b\, E(t,x;b,y;\mu,\nu^2) \Big|_{b=0}, \label{def K0(t,x;y)}\\
K_1(t,x;y;\mu,\nu^2) & \doteq  E(t,x;0,y;\mu,\nu^2) \label{def K1(t,x;y)}
\end{align} with parameter $\gamma\doteq \frac{1-\sqrt{\delta}}{2}$ and $\mathsf{F}(a,b;c; z)$  Gauss hypergeometric function.
\end{proposition}

\begin{remark} \label{Remark lower bound Hyper Gauss funct} In the next sections, we will need to estimate from below the kernel function $E$. In particular, we use the lower bound estimate
\begin{align} \label{lower bound Hyper Gauss funct}
\mathsf{F}(a,a;c;z)\geq 1 
\end{align} for  any $z\in [0,1)$ when $a\in\mathbb{R}$ and $c>0$. This estimate follows trivially from the series expansion of $\mathsf{F}(a,a;c;z)$.
\end{remark}

In the next subsection, we will prove the blow-up result by using this representation formula for an auxiliary function related to a local solution to \eqref{semilinear CP derivative type}.

\subsection{Comparison argument} \label{Subsection comparison argument}

Let us consider a local (in time) solution $u$ of the Cauchy problem \eqref{semilinear CP derivative type}. Then, we introduce a new function which depends on the time variable and only on the first space variable, by integrating with respect to the  remaining $(n-1)$ spatial variables. That is, if we denote $x=(z,w)$ with $z\in \mathbb{R}$ and $w\in \mathbb{R}^{n-1}$, then, we deal with the function
\begin{align*}
\mathcal{U}(t,z)\doteq \int_{\mathbb{R}^{n-1}} u(t,z,w) \, \mathrm{d} w \qquad \mbox{for any} \ \  t>0, z\in\mathbb{R}.
\end{align*} Of course, in the one dimensional case we may work directly with $u$ instead of $\mathcal{U}$. Hereafter, we will deal only with the case $n\geq 2$ for the sake of brevity, although one can proceed exactly in the same way for $n=1$ by working with $u$ in place of $\mathcal{U}$. Similarly, we introduce
\begin{align*}
\mathcal{U}_0(z)\doteq \int_{\mathbb{R}^{n-1}} u_0(z,w) \, \mathrm{d} w , \quad \mathcal{U}_1(z)\doteq \int_{\mathbb{R}^{n-1}} u_1(z,w) \, \mathrm{d} w\qquad \mbox{for any} \ \  z\in\mathbb{R}.
\end{align*} Since we assume that $u_0,u_1$ are compactly supported with support contained in $B_R$, it follows that $\mathcal{U}_0,\mathcal{U}_1$ are compactly supported in $[-R,R]$. Analogously, as $\supp u(t,\cdot)\subset B_{R+t}$ for any $t>0$, due to the property of finite speed of propagation of perturbations, we have $\supp \mathcal{U}(t,\cdot)\subset [-(R+t),R+t]$  for any $t>0$. 

Therefore, $\mathcal{U}$ solves the following Cauchy problem
\begin{align*} 
\begin{cases}
\mathcal{U}_{tt}- \mathcal{U}_{zz} +\frac{\mu}{1+t}\mathcal{U}_t+\frac{\nu^2}{(1+t)^2}\mathcal{U}=\int_{\mathbb{R}^{n-1}} |\partial_t u(t,z,w)|^p \, \mathrm{d} w , &  z\in \mathbb{R} , \ t>0,\\
\mathcal{U}(0,z)= \varepsilon\, \mathcal{U}_0(z), & z\in \mathbb{R} , \\ \mathcal{U}_t(0,z)= \varepsilon\, \mathcal{U}_1(z) , & z\in \mathbb{R}.
\end{cases}
\end{align*} By Proposition \ref{Prop representation formula 1d case} we know an explicit representation for $\mathcal{U}$. Also,
\begin{align*}
\mathcal{U}(t,z) &= \frac{\varepsilon}{2}(1+t)^{-\frac{\mu}{2}}\big(\mathcal{U}_0(z+t)+\mathcal{U}_0(z-t)\big)\\ & \quad +\frac{\varepsilon}{2^{\sqrt{\delta}}}\int_{z-t}^{z+t} \Big[ \mathcal{U}_0(y) K_0(t,z;y;\mu,\nu^2) + \big(\mathcal{U}_1(y)+\mu \, \mathcal{U}_0(y)\big) K_1(t,z;y;\mu,\nu^2)\Big]\, \mathrm{d}y  \\ & \quad +\frac{1}{2^{\sqrt{\delta}}} \int_0^t \int_{z-t+b}^{z+t-b}\int_{\mathbb{R}^{n-1}} |\partial_t u(b,y,w)|^p \, \mathrm{d} w \, E(t,z;b,y;\mu,\nu^2) \, \mathrm{d}y\, \mathrm{d}b,
\end{align*} where the kernel functions $E,K_0,K_1$ are defined by \eqref{def E(t,x;b,y)}, \eqref{def K0(t,x;y)} and \eqref{def K1(t,x;y)}, respectively. 

Due to the sign assumption for $u_0$ it follows that $\mathcal{U}_0$ is a nonnegative function. Consequently, from the last equality we get
\begin{align*}
\mathcal{U}(t,z) &\geq \frac{\varepsilon}{2^{\sqrt{\delta}}}\int_{z-t}^{z+t} \mathcal{U}_0(y) \Big[  K_0(t,z;y;\mu,\nu^2) + \mu K_1(t,z;y;\mu,\nu^2)\Big]\,  \mathrm{d}y  \\ 
 & \quad + \frac{\varepsilon}{2^{\sqrt{\delta}}}\int_{z-t}^{z+t} \mathcal{U}_1(y) K_1(t,z;y;\mu,\nu^2) \, \mathrm{d}y  \\
& \quad +\frac{1}{2^{\sqrt{\delta}}} \int_0^t \int_{z-t+b}^{z+t-b}\int_{\mathbb{R}^{n-1}} |\partial_t u(b,y,w)|^p \, \mathrm{d} w \, E(t,z;b,y;\mu,\nu^2) \, \mathrm{d}y\, \mathrm{d}b\doteq \varepsilon J +I.
\end{align*} 
 Let us estimate from below the two addends in the  last inequality for $\mathcal{U}(t,z)$, which are denoted by $J$ and $I$. According to Remark \ref{Remark lower bound Hyper Gauss funct}, for $\mu,\nu^2$ such that $\delta\geq 0$ the hypergeometric function that appears in the kernel $K_1$ is estimated from below by a constant for $|z-y|\leq t$. Hence, using  $$4(t+1)\leq (t+2)^2 -(y-z)^2 \leq (t+2)^2$$ for $y\in [z-t,z+t]$, we obtain
\begin{align*} 
K_1(t,z;y;\mu,\nu^2) &\gtrsim (1+t)^{-\frac{\mu}{2}+\gamma}  \left((t+2)^2-(y-z)^2\right)^{-\gamma} \notag
\gtrsim 
\begin{cases} 
(1+t)^{-\frac{\mu}{2}-\gamma} & \mbox{if} \ \delta\in[0,1), \\  (1+t)^{-\frac{\mu}{2}} & \mbox{if} \ \delta\geq 1,
\end{cases}
\end{align*}  where in the second inequality we estimated the factor containing $y$ with its minimum on $[z-t,z+t]$, that is,
\begin{align}
K_1(t,z;y;\mu,\nu^2) &\gtrsim  (1+t)^{-\frac{\sigma}{2}} \label{est below K1}
\end{align} with $\sigma$ defined by \eqref{def shift sigma}.
Elementary computations lead to 
\begin{align*}
\partial_b E(t,z; & 0,y;\mu,\nu^2) \\  & = (1+t)^{-\frac{\mu}{2}+\gamma} ((t+2)^2-(y-z)^2)^{-\gamma} \\ 
& \quad \times \Big[ 4\gamma^2 (1+t)((y-z)^2-t(t+2)) ((t+2)^2-(y-z)^2)^{-2} \,\mathsf{F}(\gamma+1,\gamma+1;2;\zeta) \\
& \qquad \quad + (\tfrac{\mu}{2}+\gamma)\, \mathsf{F}(\gamma,\gamma; 1;\zeta) -2\gamma  ((t+2)^2-(y-z)^2)^{-1}(t+2)  \,\mathsf{F}(\gamma,\gamma; 1;\zeta) \Big],
\end{align*} where $\zeta = \zeta (t,z;y)\doteq  \frac{t^2-(y-z)^2}{(t+2)^2-(y-z)^2}$. Therefore, for $y\in [z-t,z+t]$ and $\delta\geq 1$, since $\gamma\leq 0$,  we have
\begin{align*}
  K_0(t,z; & y;\mu,\nu^2)+\mu K_1(t,z;y;\mu,\nu^2) \notag \\&  = \mu  E(t,z;  0,y;\mu,\nu^2)   -\partial_b E(t,z;0,y;\mu,\nu^2)  \notag\\ & \geq  (1+t)^{-\frac{\mu}{2}+\gamma} ((t+2)^2-(y-z)^2)^{-\gamma} \Big[
 (\tfrac{\mu}{2}-\gamma)\, +  \tfrac{2\gamma(t+2)}{(t+2)^2-(y-z)^2}  \Big]  \mathsf{F}(\gamma,\gamma; 1;\zeta) \notag
\\ & \geq  (1+t)^{-\frac{\mu}{2}+\gamma} ((t+2)^2-(y-z)^2)^{-\gamma} \Big[
 (\tfrac{\mu}{2}-\gamma)\, +  \tfrac{\gamma(t+2)}{2(t+1)}  \Big]  \mathsf{F}(\gamma,\gamma; 1;\zeta), \notag
 \end{align*} where in the last inequality we estimated the second addend in brackets by its minimum. Also, for $y\in [z-t,z+t]$ and $\delta\geq 1$ we get
 \begin{align}
  K_0(t,z; & y;\mu,\nu^2)+\mu K_1(t,z;y;\mu,\nu^2) \notag \\
  &  \geq  (1+t)^{-\frac{\mu}{2}+\gamma} ((t+2)^2-(y-z)^2)^{-\gamma} \Big[
 (\tfrac{\mu}{2}-\gamma)\, +  \tfrac{\gamma}{2}\Big(2-\tfrac{t}{t+1}\Big)  \Big]  \mathsf{F}(\gamma,\gamma; 1;\zeta) \notag
\\ & \geq \tfrac{\mu}{2} (1+t)^{-\frac{\mu}{2}+\gamma} ((t+2)^2-(y-z)^2)^{-\gamma}  \mathsf{F}(\gamma,\gamma; 1;\zeta) \notag
\\ & \gtrsim (1+t)^{-\frac{\mu}{2}}= (1+t)^{-\frac{\sigma}{2}}, \label{est below K0+muK1}
\end{align} where in the second step we used $\gamma\leq 0$ and in the last step we used the same estimate from below as in \eqref{est below K1} and \eqref{lower bound Hyper Gauss funct}. 

Let us remark that $[-R,R]\subset [z-t,z+t]$ if and only if $z\in [-t+R,t-R]$. Thus, we found for $z\in [-t+R,t-R]$ (in the case $\delta\geq 1$) the estimate
\begin{align}
J=  \frac{1}{2^{\sqrt{\delta}}}\int_{z-t}^{z+t} & \Big[ \mathcal{U}_0(y) K_0(t,z;y;\mu,\nu^2) + \big(\mathcal{U}_1(y)+\mu \, \mathcal{U}_0(y)\big) K_1(t,z;y;\mu,\nu^2)\Big]\, \mathrm{d}y 
 \notag\\ & \quad \gtrsim  (1+t)^{-\frac{\sigma}{2}} \int_{z-t}^{z+t}  \big(\mathcal{U}_1(y)+ \mathcal{U}_0(y)\big)\, \mathrm{d}y \notag \\ & \quad = (1+t)^{-\frac{\sigma}{2}} \int_{\mathbb{R}}  \big(\mathcal{U}_1(y)+ \mathcal{U}_0(y)\big)\, \mathrm{d}y =  \|u_0+u_1\|_{L^1(\mathbb{R}^n)}  (1+t)^{-\frac{\sigma}{2}}. \label{lower bound estimate J int}
\end{align}
In the case $\delta\in[0,1)$, we may not prove the estimate in \eqref{est below K0+muK1} for the term $K_0+\mu K_1$ as in the previous case, due to the fact that $\gamma$ is positive. Nonetheless, \eqref{est below K1} still holds for $K_1$. Hence, assuming $u_0=0$ in the latter case, we get once again the lower bound estimate for $J$ in \eqref{lower bound estimate J int}.

Next we estimate the term $I$. Since $\supp u(t,\cdot) \subset B_{R+t}$  implies $$\supp \partial_t u(t,z,\cdot)\subset \big\{w\in \mathbb{R}^{n-1}: |w|\leq \left((R+t)^2-z^2\right)^{1/2} \big\} \qquad \mbox{for any} \ \ t>0, z\in \mathbb{R}$$ by H\"older's inequality we get
\begin{align*}
  |\mathcal{U}_t(b,y)| & =\Big| \int_{\mathbb{R}^{n-1}} \partial_t u(b, y,w) \, \mathrm{d}w\Big| \\ &  \leq \bigg(\int_{\mathbb{R}^{n-1}} |\partial_t u(b, y,w)|^p \, \mathrm{d}w\bigg)^{\frac{1}{p}} \left(\meas\big(\supp \partial_t u(b,y,\cdot)\big)\right)^{1-\frac{1}{p}} \\ & \lesssim \left((R+b)^2-y^2\right)^{\frac{n-1}{2}\left(1-\frac{1}{p}\right)}\bigg(\int_{\mathbb{R}^{n-1}} |\partial_t u(b, y,w)|^p \, \mathrm{d}w\bigg)^{\frac{1}{p}} .
\end{align*} Hence,
\begin{align*}
\int_{\mathbb{R}^{n-1}} |\partial_t u(b, y,w)|^p \, \mathrm{d}w  & \gtrsim  \left((R+b)^2-y^2\right)^{-\frac{n-1}{2}(p-1)} |\mathcal{U}_t(b,y)|^p,
\end{align*}  which implies in turn
\begin{align*}
I &  \gtrsim  \int_0^t \int_{z-t+b}^{z+t-b}\left((R+b)^2-y^2\right)^{-\frac{n-1}{2}(p-1)} |\mathcal{U}_t(b,y)|^p  E(t,z;b,y;\mu,\nu^2) \, \mathrm{d}y\, \mathrm{d}b \\ 
&=  \int_{z-t}^{z+t} \int_0^{t-|y-z|} \left((R+b)^2-y^2\right)^{-\frac{n-1}{2}(p-1)} |\mathcal{U}_t(b,y)|^p  E(t,z;b,y;\mu,\nu^2)  \, \mathrm{d}b \, \mathrm{d}y ,
\end{align*} where we used Fubini's theorem in the last equality. We work now on the characteristic $t-z=R$. Then, shrinking the domain of integration, we have
\begin{align*}
I &  \gtrsim \int_{R}^{z} \int_{y-R}^{y+R} \left((R+b)^2-y^2\right)^{-\frac{n-1}{2}(p-1)} |\mathcal{U}_t(b,y)|^p  E(t,z;b,y;\mu,\nu^2)  \, \mathrm{d}b \, \mathrm{d}y \\
& \gtrsim \int_{R}^{z} \left(R+y\right)^{-\frac{n-1}{2}(p-1)} \int_{y-R}^{y+R}|\mathcal{U}_t(b,y)|^p  E(t,z;b,y;\mu,\nu^2)  \, \mathrm{d}b \, \mathrm{d}y.
\end{align*} Note that the unexpressed multiplicative constant in the previous chain of inequalities depends on $R$. Now we estimate from below the kernel function $E$. 
 Consequently, using again \eqref{lower bound Hyper Gauss funct}, for $b\in [0,t], y\in [z-(t-b),z+t-b]$ we may estimate
\begin{align*}
E(t,z;b,y;\mu,\nu^2) & \gtrsim   (1+t)^{-\frac{\mu}{2}+\gamma} (1+b)^{\frac{\mu}{2}+\gamma} \left((t+b+2)^2-(y-z)^2\right)^{-\gamma} \\ &  \gtrsim   (1+t)^{-\frac{\sigma}{2}} (1+b)^{\frac{\sigma}{2}}. 
\end{align*} Also, on the characteristic $t-z=R$ we get
\begin{align*}
I &  \gtrsim  (1+t)^{-\frac{\sigma}{2}}  \int_{R}^{z} \left(R+y\right)^{-\frac{n-1}{2}(p-1)} \int_{y-R}^{y+R}|\mathcal{U}_t(b,y)|^p  (1+b)^{\frac{\sigma}{2}} \, \mathrm{d}b \, \mathrm{d}y \\
&  \gtrsim  (1+t)^{-\frac{\sigma}{2}}  \int_{R}^{z} \left(R+y\right)^{-\frac{n-1}{2}(p-1)+\frac{\sigma}{2}} \int_{y-R}^{y+R}|\mathcal{U}_t(b,y)|^p  \left(\frac{1+b}{R+y}\right)^{\frac{\sigma}{2}}  \mathrm{d}b \, \mathrm{d}y.
\end{align*} We notice that the quotient in the $b$-integral in the last line is bounded from below on the domain of integration by a positive constant, that depends on $R$. Clearly, we can assume without loss of generality $R>1$. We take $y\in [R,z]$ and $b\in [y-R,y+R]$. Then,
\begin{align*}
\left(\frac{1+b}{R+y}\right)^{\frac{\sigma}{2}}   \geq \left(\frac{1+y-R}{R+y}\right)^{\frac{\sigma}{2}}  \geq C_R.
\end{align*} The last inequality can be proved by splitting the cases $y\in [R,3R-2]$ and  $y\geq 3R-2$, as follows:
\begin{align*}
\left(\frac{1+y-R}{R+y}\right)^{\frac{\sigma}{2}} & \geq (R+y)^{-\frac{\sigma}{2}} \geq (2(2R-1))^{-\frac{\sigma}{2}} \qquad \mbox{for} \ \ y\in [R,3R-2], \\
\left(\frac{1+y-R}{R+y}\right)^{\frac{\sigma}{2}} & \geq 2^{-\frac{\sigma}{2}} \qquad \qquad \qquad \qquad  \qquad \qquad \ \  \mbox{for} \ \ y\geq 3R-2.
\end{align*} Therefore, by using Jensen's inequality and the fundamental theorem of calculus we arrive at
\begin{align*}
I  & \gtrsim  (1+t)^{-\frac{\sigma}{2}}  \int_{R}^{z} \left(R+y\right)^{-\frac{n-1}{2}(p-1)+\frac{\sigma}{2}} \int_{y-R}^{y+R}|\mathcal{U}_t(b,y)|^p   \mathrm{d}b \, \mathrm{d}y \\
& \gtrsim  (1+t)^{-\frac{\sigma}{2}}  \int_{R}^{z} \left(R+y\right)^{-\frac{n-1}{2}(p-1)+\frac{\sigma}{2}} \bigg|\int_{y-R}^{y+R}\mathcal{U}_t(b,y)   \mathrm{d}b \bigg|^p  \mathrm{d}y  \\
& = (1+t)^{-\frac{\sigma}{2}}  \int_{R}^{z} \left(R+y\right)^{-\frac{n-1}{2}(p-1)+\frac{\sigma}{2}}|\mathcal{U}(y+R,y)  |^p \, \mathrm{d}y \\
& = (1+t)^{-\frac{\sigma}{2}}  \int_{R}^{z} \left(R+y\right)^{-\frac{n+\sigma-1}{2}(p-1)}\left|(R+y)^{\frac{\sigma}{2}}\mathcal{U}(y+R,y) \right |^p \mathrm{d}y,
\end{align*} where in the third step we used $\mathcal{U}(y-R,y)=0$. Combining the lower bound estimates for $J$ and $I$, on the characteristic $t-z=R$ and for $t\geq 2R$ we found
\begin{align*}
(R+z)^{\frac{\sigma}{2}}\mathcal{U}(z+R,z) \gtrsim \varepsilon \, \| u_0+u_1\|_{L^1(\mathbb{R}^n)} \! + \! \int_{R}^{z} \! \! \left(R+y\right)^{-\frac{n+\sigma-1}{2}(p-1)}\! \left|(R+y)^{\frac{\sigma}{2}}\mathcal{U}\!(y+R,y) \right |^p \!  \mathrm{d}y.
\end{align*} If we introduce the function $U(z)\doteq (R+z)^{\frac{\sigma}{2}} \mathcal{U}(z+R,z)$ and we denote by $C$ the unexpressed multiplicative constant in the last inequality, we may rewrite  
\begin{align}\label{fundamental inequality for U}
U(z)\geq M  \varepsilon +  C \int_{R}^{z} \left(R+y\right)^{-\frac{n+\sigma-1}{2}(p-1)} |U(y) |^p  \mathrm{d}y \qquad \mbox{for any} \ \ z\geq R,
\end{align} where $M \doteq C\| u_0+u_1\|_{L^1(\mathbb{R}^n)} $.
Let us introduce the function 
\begin{align*}
G(z)\doteq   M  \varepsilon +  C \int_{R}^{z} \left(R+y\right)^{-\frac{n+\sigma-1}{2}(p-1)} |U(y) |^p  \mathrm{d}y \qquad \mbox{for any} \ \ z\geq R.
\end{align*} Clearly, by \eqref{fundamental inequality for U} we obtain $U\geq G$. Moreover, $G$ solves the differential inequality
\begin{align*}
G'(z) & = C  \left(R+z\right)^{-\frac{n+\sigma-1}{2}(p-1)} |U(z) |^p \\
& \geq C  \left(R+z\right)^{-\frac{n+\sigma-1}{2}(p-1)} (G(z))^p. 
\end{align*} As $G$ is a positive function, then, separation of variables leads to 
\begin{align*}
\tfrac{(M\varepsilon)^{1-p}-G(z)^{1-p} }{p-1}\geq  \tfrac{C}{1-\tfrac{n+\sigma-1}{2}(p-1)} \Big((R+z)^{-\frac{n+\sigma-1}{2}(p-1)+1}-(2R)^{-\frac{n+\sigma-1}{2}(p-1)+1}\Big) 
\end{align*} in the subcritical case $p\in (1,p_{\Gla}(n+\sigma))$ and 
\begin{align*}
\tfrac{(M\varepsilon)^{1-p}-G(z)^{1-p} }{p-1}\geq C\log \left(\frac{R+z}{2R}\right) 
\end{align*} if $p=p_{\Gla}(n+\sigma) $.
 In the subcritical case $p\in (1,p_{\Gla}(n+\sigma))$, choosing $\varepsilon\in (0,\varepsilon_0]$ sufficiently small with $\varepsilon_0=\varepsilon_0(n,p,\mu,\nu^2,u_0,u_1,R)$, we get 
\begin{align*}
G(z) & \geq \left[(M\varepsilon)^{1-p}+\frac{C}{\frac{1}{p-1}-\frac{n+\sigma-1}{2}} \Big((2R)^{-\frac{n+\sigma-1}{2}(p-1)+1}-(R+z)^{-\frac{n+\sigma-1}{2}(p-1)+1}\Big)\right]^{-\frac{1}{p-1}}  \\
&\geq \left[2(M\varepsilon)^{1-p}-\frac{C}{\frac{1}{p-1}-\frac{n+\sigma-1}{2}} (R+z)^{-\frac{n+\sigma-1}{2}(p-1)+1}\right]^{-\frac{1}{p-1}}.
\end{align*} From this last estimate we see that for $t=R+z\simeq \varepsilon^{-\left(\frac{1}{p-1}-\frac{n+\sigma-1}{2}\right)}$ the  lower bound for $F$  blows up.
 Then, $G$ (and $U$ in turn) blows up in finite time and the upper bound for the lifespan $$T(\varepsilon)\lesssim  \varepsilon^{-\left(\frac{1}{p-1}-\frac{n+\sigma-1}{2}\right)}$$ is fulfilled in the subcritical case. Analogously, in the critical case $p=p_{\Gla}(n+\sigma)$ we have that
 \begin{align*}
 G(z) & \geq  \left[(M\varepsilon)^{1-p}-C(p-1)\log \left(\tfrac{R+z}{2R}\right)\right]^{-\frac{1}{p-1}}
 \end{align*} implies the blow-up in finite time of $U(z)$  and the lifespan estimate $$T(\varepsilon)\lesssim \exp\left(\widetilde{C}\varepsilon^{-(p-1)}\right).$$ So, the proof of Theorem \ref{Thm blow-up derivative type nonlinearity} is complete.
 
 \section{Proof of Theorem \ref{Thm blow-up WCS der type nonlinearities}}
 
 In this section we prove the blow-up result for the weakly coupled system \eqref{semilinear WCS derivative type}.
 The section is organized as follows: in Subsection \ref{Subsection iteration frame} we introduce two suitable functions which are related to the components of a local in time solution of \eqref{semilinear WCS derivative type} and we derive the corresponding iteration frame by using the same ideas from Subsection \ref{Subsection comparison argument}; then, in order to prove Theorem \ref{Thm blow-up WCS der type nonlinearities} we apply an iteration argument both in the subcritical case (Subsection \ref{Subsection iteration argument subcrit case})
 and in the critical case (Subsection \ref{Subsection iteration argument crit case}). In particular, in the critical case we employ the so-called \emph{slicing method} in order to deal with  logarithmic factors. For further details on the slicing method see \cite{AKT00}, where this method was introduced for the first time or \cite{TakWak11,TakWak14,WakYor18,PalTak19,PalTak19dt,PalTak19mix}  where the slicing method is used in critical cases in order to manage factors of logarithmic type.

 \subsection{Iteration frame} \label{Subsection iteration frame}
 
Let $(u,v)$ be a local in time solution to \eqref{semilinear WCS derivative type}. If we denote $x=(z,w)$ with $z\in\mathbb{R}$ and $w\in \mathbb{R}^{n-1}$ as in Subsection \ref{Subsection comparison argument}, then, we may introduce the functions
\begin{align*}
\mathcal{U}(t,z)\doteq \int_{\mathbb{R}^{n-1}} u(t,z,w) \, \mathrm{d} w ,\quad \mathcal{V}(t,z)\doteq \int_{\mathbb{R}^{n-1}} v(t,z,w) \, \mathrm{d} w
\end{align*} for any $t>0, z\in\mathbb{R}$ in the case $n\geq 2$. Clearly, also in this case we can simply work with $u,v$ instead of $\mathcal{U},\mathcal{V}$ for $n=1$. Repeating the same steps as in the case of the single semilinear equation, we end up with the estimates
\begin{equation} \label{relation mathcal U,V pre}
\begin{split}
(R+z)^{\frac{\sigma_1}{2}}\mathcal{U}(z+R,z) & \gtrsim \varepsilon \, \| u_0+u_1\|_{L^1(\mathbb{R}^n)} \\ & \qquad +  \int_{R}^{z} \left(R+y\right)^{-\frac{n-1}{2}(p-1)+\frac{\sigma_1}{2}-\frac{\sigma_2}{2}p}\left|(R+y)^{\frac{\sigma_2}{2}}\mathcal{V}(y+R,y) \right |^p  \mathrm{d}y,\\
(R+z)^{\frac{\sigma_2}{2}}\mathcal{V}(z+R,z) & \gtrsim \varepsilon \, \| v_0+v_1\|_{L^1(\mathbb{R}^n)} \\ & \qquad  +  \int_{R}^{z} \left(R+y\right)^{-\frac{n-1}{2}(q-1)+\frac{\sigma_2}{2}-\frac{\sigma_1}{2}q}\left|(R+y)^{\frac{\sigma_1}{2}}\mathcal{U}(y+R,y) \right |^q  \mathrm{d}y 
\end{split}
\end{equation}
on the characteristic $t=z+R$ for $z\geq R$. Let us point out that the assumptions on the Cauchy data in the statement of Theorem \ref{Thm blow-up WCS der type nonlinearities} allow us to proceed exactly as the proof of Theorem \ref{Thm blow-up derivative type nonlinearity} when we estimate from below the terms which are related to the solution of the corresponding linear homogeneous problem. We define the functions $U(z)\doteq (R+z)^{\frac{\sigma_1}{2}}\mathcal{U}(z+R,z)$ and $V(z)\doteq (R+z)^{\frac{\sigma_2}{2}}\mathcal{V}(z+R,z)$. Hence, denoting by $C$ and $K$ the unexpressed multiplicative constants in \eqref{relation mathcal U,V pre}, we obtain the iteration frame 
\begin{align}
U(z) \geq M\varepsilon +C\int_R^z  \left(R+y\right)^{-\frac{n-1}{2}(p-1)+\frac{\sigma_1}{2}-\frac{\sigma_2}{2}p}|V(y)|^p \, \mathrm{d}y, \label{iteration frame lower bound U} \\ 
V(z) \geq N\varepsilon +K\int_R^z  \left(R+y\right)^{-\frac{n-1}{2}(q-1)+\frac{\sigma_2}{2}-\frac{\sigma_1}{2}q}|U(y)|^q \, \mathrm{d}y \label{iteration frame lower bound V} 
\end{align} for any $z\geq R$, where $M\doteq C\| u_0+u_1\|_{L^1(\mathbb{R}^n)}$ and  $N\doteq K\| v_0+v_1\|_{L^1(\mathbb{R}^n)}$.
Note that \eqref{iteration frame lower bound U} and \eqref{iteration frame lower bound V} provide not only the iteration frame for the pair $(U,V)$, but also the base step of the inductive argument. Indeed, in the base case we will simply estimate $U,V$ from below by the two quantities $M\varepsilon,N\varepsilon$, respectively.
 
 \subsection{Iteration argument: subcritical case} \label{Subsection iteration argument subcrit case}
 
 In this section we prove that a local in time solution $(u,v)$ to \eqref{semilinear WCS derivative type} blows up in finite time in the subcritical case 
 \begin{align*}
 \Omega(n,\sigma_1,\sigma_2,p,q)& =\max\left\{\Lambda(n+\sigma_1,p,q), \Lambda(n+\sigma_2,q,p)\right\}\\ &=\max\left\{\frac{p+1}{pq-1}-\frac{n+\sigma_1-1}{2}, \frac{q+1}{pq-1}-\frac{n+\sigma_2-1}{2}\right\}>0,
 \end{align*}
 of course, provided that $u_0,u_1,v_0,v_1$ satisfy the assumptions of Theorem \ref{Thm blow-up WCS der type nonlinearities}.
 
Let us assume that $\Omega(n,\sigma_1,\sigma_2,p,q)=\Lambda(n+\sigma_1,p,q)$. 
First we prove the sequence of lower bound estimates for $U$
\begin{align}\label{lower bound U j subcritical}
U(z)\geq C_j (R+z)^{-\alpha_j}(z-R)^{\beta_j} \qquad \mbox{for any} \ \ z\geq R,
\end{align} where $\{\alpha_j\}_{j\in \mathbb{N}}$, $\{\beta_j\}_{j\in \mathbb{N}}$ and $\{C_j\}_{j\in \mathbb{N}}$ are sequences of nonnegative real numbers that we will determine afterwards via an inductive procedure. Clearly, from \eqref{iteration frame lower bound U} we see that \eqref{lower bound U j subcritical} is true for $j=0$, provided that $\alpha_0\doteq 0$, $\beta_0\doteq 0$ and $C_0\doteq M\varepsilon$. We prove now the inductive step. We assume that \eqref{lower bound U j subcritical} is satisfied for $j\geq 0$.
Plugging \eqref{lower bound U j subcritical} in \eqref{iteration frame lower bound V}, we get 
\begin{align*}
V(z) & \geq K\int_R^z  \left(R+y\right)^{-\frac{n-1}{2}(q-1)+\frac{\sigma_2}{2}-\frac{\sigma_1}{2}q}|U(y)|^q \, \mathrm{d}y \\
& \geq K C_j^q \int_R^z  \left(R+y\right)^{-\frac{n-1}{2}(q-1)+\frac{\sigma_2}{2}-\frac{\sigma_1}{2}q-q \alpha_j}(y-R)^{q\beta_j} \, \mathrm{d}y \\
& \geq K C_j^q \left(R+z\right)^{-\frac{n-1}{2}(q-1)-\frac{\sigma_1}{2}q-q \alpha_j} \int_R^z  (y-R)^{q\beta_j+\frac{\sigma_2}{2}} \, \mathrm{d}y \\
& = K C_j^q \big(q\beta_j+\tfrac{\sigma_2}{2}+1\big)^{-1}\left(R+z\right)^{-\frac{n-1}{2}(q-1)-\frac{\sigma_1}{2}q-q \alpha_j}   (z-R)^{q\beta_j+\frac{\sigma_2}{2}+1} 
\end{align*}for $z\geq R$. Combining the above lower bound for $V$ and \eqref{iteration frame lower bound U}, we arrive at
\begin{align*}
U (z) & \geq C\int_R^z  \left(R+y\right)^{-\frac{n-1}{2}(p-1)+\frac{\sigma_1}{2}-\frac{\sigma_2}{2}p}|V(y)|^p \, \mathrm{d}y \\
& \geq  \frac{C K^p C_j^{pq}}{ \big(q\beta_j+\tfrac{\sigma_2}{2}+1\big)^{p}}\int_R^z \left(R+y\right)^{-\frac{n-1}{2}(pq-1)-\frac{\sigma_1}{2}pq-\frac{\sigma_2}{2}p-pq \alpha_j}  (y-R)^{pq\beta_j+\frac{\sigma_2}{2}p+p+\frac{\sigma_1}{2}} \, \mathrm{d}y \\
& \geq  \frac{C K^p C_j^{pq}}{ \big(q\beta_j+\tfrac{\sigma_2}{2}+1\big)^{p}} \left(R+z\right)^{-\frac{n-1}{2}(pq-1)-\frac{\sigma_1}{2}pq-\frac{\sigma_2}{2}p-pq \alpha_j} \int_R^z   (y-R)^{pq\beta_j+\frac{\sigma_2}{2}p+\frac{\sigma_1}{2}+p} \, \mathrm{d}y \\
& =  \frac{C K^p C_j^{pq} }{\big(q\beta_j+\tfrac{\sigma_2}{2}+1\big)^{p} \big(pq\beta_j+\tfrac{\sigma_2}{2}p+\tfrac{\sigma_1}{2}+p+1\big)}   \left(R+z\right)^{-\frac{n-1}{2}(pq-1)-\frac{\sigma_1}{2}pq-\frac{\sigma_2}{2}p-pq \alpha_j} \\ & \qquad \qquad \qquad \qquad \qquad\qquad \qquad \qquad \qquad \quad \times (z-R)^{pq\beta_j+\frac{\sigma_2}{2}p+\frac{\sigma_1}{2}+p+1} 
\end{align*} for $z\geq R$. So, we proved \eqref{lower bound U j subcritical}  for $j+1$, provided that 
\begin{align}
\alpha_{j+1} &\doteq \tfrac{n+\sigma_1-1}{2}(pq-1)+\tfrac{\sigma_2}{2}p+\tfrac{\sigma_1}{2}+pq\, \alpha_j , \label{def alpha j+1}\\
\beta_{j+1} &\doteq \tfrac{\sigma_2}{2}p+\tfrac{\sigma_1}{2}+p+1+pq\, \beta_j , \label{def beta j+1}\\
C_{j+1}& \doteq CK^p \big( \tfrac{\sigma_2}{2}+1 + q\, \beta_j \big)^{-p} \big( \tfrac{\sigma_2}{2}p+\tfrac{\sigma_1}{2}+p+1+pq\, \beta_j \big)^{-1}C_j^{pq}. \label{def C j+1}
\end{align} Next we derive the explicit expressions for $\alpha_j$ and $\beta_j$. Applying \eqref{def alpha j+1} iteratively, we get
\begin{align}\label{expression alpha j}
\alpha_j = A+ pq\, \alpha_{j-1} = A(1+pq)+ (pq)^2 \alpha_{j-2}= A\sum_{k=0}^{j-1} (pq)^k+(pq)^j \alpha_0 =A \tfrac{(pq)^j-1}{pq-1} ,
\end{align} where $A\doteq \tfrac{n+\sigma_1-1}{2}(pq-1)+\tfrac{\sigma_2}{2}p+\tfrac{\sigma_1}{2}$ and we used $\alpha_0=0$. Similarly, from \eqref{def beta j+1} we find
\begin{align} \label{expression beta j}
\beta_j = B \tfrac{(pq)^j-1}{pq-1} ,
\end{align} where $B\doteq \tfrac{\sigma_2}{2}p+\tfrac{\sigma_1}{2}+p+1$.  Now we use the explicit expression of $\beta_j$ to get a lower bound estimate for $C_j$. 
Since $\beta_j\leq \frac{B}{pq-1}(pq)^j$, by \eqref{def beta j+1} we get 
\begin{align*}
C_j & =CK^p \big( \tfrac{\sigma_2}{2}+1 + q\, \beta_{j-1} \big)^{-p} \big( \tfrac{\sigma_2}{2}p+\tfrac{\sigma_1}{2}+p+1+pq\, \beta_{j-1} \big)^{-1}C_{j-1}^{pq} \\ & \geq CK^p  \big( \tfrac{\sigma_2}{2}p+\tfrac{\sigma_1}{2}+p+1+pq\, \beta_{j-1} \big)^{-(p+1)}C_{j-1}^{pq} \\ & =  CK^p \beta_j^{-(p+1)} C_{j-1}^{pq} \geq \underbrace{ CK^p \left(\tfrac{B}{pq-1}\right)^{-(p+1)}}_{\doteq \, \widetilde{C}} (pq)^{-(p+1)j}   C_{j-1}^{pq} =\widetilde{C} (pq)^{-(p+1)j}   C_{j-1}^{pq}.
\end{align*} Applying  the logarithmic function to both sides of the last inequality and using iteratively the resulting inequality, we get 
\begin{align*}
\log C_j & \geq pq \log C_{j-1} -j \log (pq)^{p+1}+\log \widetilde{C} \\
& \geq (pq)^2 \log C_{j-2} -(j+(j-1)(pq)) \log (pq)^{p+1}+(1+pq)\log \widetilde{C} \\
& \geq \cdots \geq  (pq)^j \log C_{0} -\left(\sum_{k=0}^{j-1}(j-k)(pq)^k\right) \log (pq)^{p+1}+\left(\sum_{k=0}^{j-1}(pq)^k\right)\log \widetilde{C}.
\end{align*} Using the formula
\begin{align*}
\sum_{k=0}^{j-1}(j-k)(pq)^k = \frac{1}{pq-1}\left(\frac{(pq)^{j+1}-1}{pq-1}-(j+1)\right),
\end{align*} which can be proved with an inductive argument, we have
\begin{align*}
\log C_j & \geq  (pq)^j \log C_{0} -\frac{1}{pq-1}\left(\frac{(pq)^{j+1}-1}{pq-1}-(j+1)\right)\log (pq)^{p+1}+\frac{(pq)^{j}-1}{pq-1}\log \widetilde{C} \\
&= (pq)^j\left( \log C_{0} -\frac{pq}{(pq-1)^2}\log (pq)^{p+1}+\frac{\log \widetilde{C}}{pq-1}\right) +(j+1)\frac{\log (pq)^{p+1}}{pq-1} \\ & \qquad +\frac{\log (pq)^{p+1}}{(pq-1)^2}-\frac{\log \widetilde{C}}{pq-1}.
\end{align*} For $j\geq j_0\doteq \max\left\{0,\frac{\log\widetilde{C}}{\log(pq)^{p+1}}-\frac{pq}{pq-1}\right\}$ from the last inequality we get 
\begin{align}\label{lower bound Cj}
\log C_j \geq (pq)^j \log(\widehat{C} \varepsilon),
\end{align} where $\widehat{C}\doteq M (pq)^{-\frac{(pq)(p+1)}{(pq-1)^2}}\widetilde{C}^{\frac{1}{pq-1}}$. Finally, we combine \eqref{lower bound U j subcritical}, \eqref{expression alpha j}, \eqref{expression beta j} and \eqref{lower bound Cj} and it results
\begin{align}
U(z) & \geq C_j (R+z)^{-A \frac{(pq)^j-1}{pq-1}} (z-R)^{B \frac{(pq)^j-1}{pq-1}}\notag \\ &\geq \exp \left((pq)^j \log(\widehat{C} \varepsilon)\right) (R+z)^{-A \frac{(pq)^j-1}{pq-1}} (z-R)^{B \frac{(pq)^j-1}{pq-1}} \notag \\
&=\exp \left((pq)^j \! \left( \log(\widehat{C} \varepsilon)-\tfrac{A}{pq-1}\log(R+z)+\tfrac{B}{pq-1}\log(z-R)\right)\! \right) \! (R+z)^{ \frac{A}{pq-1}}\! (z-R)^{ -\frac{B}{pq-1}} \label{lower bound U intermediate}
\end{align} for $z\geq R$ and $j\geq j_0$. If  we require $z\geq 3R$, then, it holds $2(z-R)\geq R+z$. So, from \eqref{lower bound U intermediate} we get
\begin{align*}
U(z) & \geq \exp \left((pq)^j  \log\left(2^{-\frac{B}{pq-1}}\widehat{C} \varepsilon(R+z)^{\frac{B-A}{pq-1}}\right)\right) (R+z)^{ \frac{A}{pq-1}} (z-R)^{ -\frac{B}{pq-1}} \\
& = \exp \left((pq)^j  \log\left(\bar{C} \varepsilon(R+z)^{\frac{p+1}{pq-1}-\frac{n+\sigma_1-1}{2}}\right)\right) (R+z)^{ \frac{A}{pq-1}} (z-R)^{ -\frac{B}{pq-1}}
\end{align*} for $z\geq 3R$ and $j\geq j_0$, where $\bar{C}\doteq 2^{-\frac{B}{pq-1}}\widehat{C}$. We recall that we are working for $(t,z)$ on the characteristic $t=z+R$, so we may rewrite the last inequality as
\begin{align} \label{lower bound U final}
U(z) & \geq \exp \left((pq)^j  \log\left(\bar{C} \varepsilon\,  t^{\Lambda(n+\sigma_1,p,q)}\right)\right) (R+z)^{ \frac{A}{pq-1}} (z-R)^{ -\frac{B}{pq-1}}
\end{align} for $t\geq 4R$ and $j\geq j_0$. We choose $\varepsilon_0=\varepsilon_0(n,\mu_1,\mu_2,\nu_1^2,\nu_2^2,u_0,u_1,v_0,v_1,R)>0$ such that $$(\bar{C} \varepsilon_0)^{-(\Lambda(n+\sigma_1,p,q))^{-1}}\geq 4R .$$ Then, for any $\varepsilon\in (0,\varepsilon_0]$ and $t> (\bar{C} \varepsilon)^{-(\Lambda(n+\sigma_1,p,q))^{-1}}$ we obtain
\begin{align*}
t\geq 4R \ \ \mbox{and} \ \ \bar{C} \varepsilon\,  t^{\Lambda(n+\sigma_1,p,q)}> 1
\end{align*} and, hence, letting $j\to \infty$ in \eqref{lower bound U final} the lower bound for $U(z)$ blows up. Therefore, in order to get a finite value of $U(z)$, it must hold the converse inequality for $t$. So, we have showed the upper bound for the  lifespan $$T(\varepsilon) \lesssim \varepsilon^{-(\Lambda(n+\sigma_1,p,q))^{-1}}.$$ In the case $\Omega(n,\sigma_1,\sigma_2,p,q)=\Lambda(n+\sigma_2,q,p)$ it suffices to switch the role of $U$ and $V$ in order to show the estimate $T(\varepsilon) \lesssim \varepsilon^{-(\Lambda(n+\sigma_2,q,p))^{-1}}$ in an analogous way. Also, we completed the proof of Theorem \ref{Thm blow-up WCS der type nonlinearities} in the subcritical case. In the critical case $\Omega(n,\sigma_1,\sigma_2,p,q)=0$ we need to modify our approach. As we have already announced, we will employ the slicing method in order to deal with  logarithmic factors in the sequence of lower bounds for $U$.
 
 \subsection{Iteration argument: critical case} \label{Subsection iteration argument crit case}
 
 In this subsection we prove Theorem \ref{Thm blow-up WCS der type nonlinearities} in the critical case
\begin{align*}
\Omega(n,\sigma_1,\sigma_2,p,q)& =\max\left\{\Lambda(n+\sigma_1,p,q), \Lambda(n+\sigma_2,q,p)\right\}=0.
\end{align*}
We begin with the case $\Lambda(n+\sigma_1,p,q)=0> \Lambda(n+\sigma_2,q,p)$.

Let us introduce the succession $\{\ell_j\}_{j\in \mathbb{N}}$, where $\ell_j\doteq 2-2^{-(j+1)}$. Our goal is to prove the sequence of lower bound estimate for $U$
\begin{align}\label{lower bound U j critical}
U(z)\geq D_j \left(\log\left(\frac{z}{\ell_j R}\right)\right)^{\theta_j} \qquad \mbox{for} \ \ z\geq \ell_j R, 
\end{align} where $\{D_j\}_{j\in\mathbb{N}}$, $\{\theta_j\}_{j\in\mathbb{N}}$ are suitable sequences of nonnegative real numbers that we shall determine throughout the iteration procedure. Obviously, \eqref{lower bound U j critical} is true for $j=0$ provided that $D_0\doteq M\varepsilon$ and $\theta_0\doteq 0$. Also, we proved the base case. It remains to prove the inductive step. Before starting we remark that $\{\ell_j\}_{j\in \mathbb{N}}$ is an increasing and bounded sequence. In particular, $\ell_j\geq \ell_0= \frac{3}{2}$. Consequently, for any $j\in \mathbb{N}$ and any $z\geq \ell_j R$ we may use the inequality $z\geq \frac{3}{5} (R+z)$. Let us assume that \eqref{lower bound U j critical} holds, we shall prove that \eqref{lower bound U j critical} is satisfied also for $j+1$. Combining \eqref{iteration frame lower bound V} and \eqref{lower bound U j critical}, we get 
\begin{align*}
V(z) & \geq K\int_{\ell_j R}^z  \left(R+y\right)^{-\frac{n-1}{2}(q-1)+\frac{\sigma_2}{2}-\frac{\sigma_1}{2}q}|U(y)|^q \, \mathrm{d}y \\
& \geq K D_j^q \int_{\ell_j R}^z  \left(R+y\right)^{-\frac{n-1}{2}(q-1)+\frac{\sigma_2}{2}-\frac{\sigma_1}{2}q} \left(\log\left(\tfrac{y}{\ell_j R}\right)\right)^{q \theta_j} \, \mathrm{d}y \\
& \geq K D_j^q \left(R+z\right)^{-\frac{n-1}{2}(q-1)-\frac{\sigma_1}{2}q} \int_{\ell_j R}^z  \left(R+y\right)^{
\frac{\sigma_2}{2}} \left(\log\left(\tfrac{y}{\ell_j R}\right)\right)^{q\theta_j} \, \mathrm{d}y \\
& \geq K D_j^q \left(R+z\right)^{-\frac{n-1}{2}(q-1)-\frac{\sigma_1}{2}q} \int_{\tfrac{\ell_j z}{\ell_{j+1}} }^z  \left(R+y\right)^{
\frac{\sigma_2}{2}} \left(\log\left(\tfrac{y}{\ell_j R}\right)\right)^{q \theta_j} \, \mathrm{d}y \\
& \geq K D_j^q \left(R+z\right)^{-\frac{n-1}{2}(q-1)-\frac{\sigma_1}{2}q}  \left(R+\tfrac{\ell_j z}{\ell_{j+1}}\right)^{
\frac{\sigma_2}{2}} \left(\log\left(\tfrac{z}{\ell_{j+1} R}\right)\right)^{q \theta_j} \left(1-\tfrac{\ell_{j}}{\ell_{j+1}}\right)z \\
& \geq \tfrac{3}{5} K D_j^q  \left(\tfrac{\ell_j }{\ell_{j+1}}\right)^{
\frac{\sigma_2}{2}} \left(1-\tfrac{\ell_{j}}{\ell_{j+1}}\right) \left(R+z\right)^{-\frac{n-1}{2}(q-1)-\frac{\sigma_1}{2}q+\frac{\sigma_2}{2}+1}  \left(\log\left(\tfrac{z}{\ell_{j+1} R}\right)\right)^{q \theta_j} 
\end{align*} for any $z\geq \ell_{j+1}R$. If we plug this lower bound for $V$ in \eqref{iteration frame lower bound U} and we use the critical condition $\Lambda(n+\sigma_1,p,q)=0$, it results
\begin{align*}
U& (z)  \geq C\int_{\ell_{j+1} R}^z  \left(R+y\right)^{-\frac{n-1}{2}(p-1)+\frac{\sigma_1}{2}-\frac{\sigma_2}{2}p}|V(y)|^p \, \mathrm{d}y \\
 & \geq   (\tfrac{3}{5})^p C  K^p D_j^{pq}\!  \left(\tfrac{\ell_j }{\ell_{j+1}}\right)^{
\frac{\sigma_2 p}{2}}\!\! \left(1-\tfrac{\ell_{j}}{\ell_{j+1}}\right)^p \!\! \int_{\ell_{j+1} R}^z \! \left(R+y\right)^{-\frac{n+\sigma_1-1}{2}(pq-1)+p} \left(\log\left(\tfrac{y}{\ell_{j+1} R}\right)\right)^{pq \theta_j}   \mathrm{d}y \\
 & \geq (\tfrac{3}{5})^{p+1} C  K^p D_j^{pq}  \left(\tfrac{\ell_j }{\ell_{j+1}}\right)^{
\frac{\sigma_2 p}{2}} \left(1-\tfrac{\ell_{j}}{\ell_{j+1}}\right)^p\int_{\ell_{j+1} R}^z \tfrac1y \left(\log\left(\tfrac{y}{\ell_{j+1} R}\right)\right)^{pq \theta_j}  \, \mathrm{d}y \\
& = (\tfrac{3}{5})^{p+1} C  K^p D_j^{pq}  \left(\tfrac{\ell_j }{\ell_{j+1}}\right)^{
\frac{\sigma_2 p}{2}} \left(1-\tfrac{\ell_{j}}{\ell_{j+1}}\right)^p (pq \,\theta_j+1)^{-1} \left(\log\left(\tfrac{z}{\ell_{j+1} R}\right)\right)^{pq \theta_j+1} 
\end{align*} for any $z\geq \ell_{j+1}R$. Since $1-\frac{\ell_j}{\ell_{j+1}}\geq 2^{-(j+3)}$ and $2 \ell_j >\ell_{j+1}$ for any $j\in \mathbb{N}$, we find
\begin{align*}
U(z) & \geq  C  K^p 2^{-jp-((4+\frac{\sigma_2}{2})p+1)} D_j^{pq}   (pq \, \theta_j+1)^{-1} \left(\log\left(\tfrac{z}{\ell_{j+1} R}\right)\right)^{pq \theta_j+1} 
\end{align*} for any $z\geq \ell_{j+1}R$. Also, we proved \eqref{lower bound U j critical} for $j+1$, provided that 
\begin{align}
\theta_{j+1} & \doteq 1+ pq \theta_j,\label{def sigma j+1} \\
D_{j+1} & \doteq C  K^p 2^{-jp-((4+\frac{\sigma_2}{2})p+1)}  (pq \, \theta_j+1)^{-1} D_j^{pq} . \label{def D j+1}
\end{align} Applying recursively \eqref{def sigma j+1}, we obtain
\begin{align}
\theta_j = 1+  pq \, \theta_{j-1} = 1+  pq+ (pq)^2 \, \theta_{j-2} = \cdots = \sum_{k=0}^{j-1}(pq)^k +\theta_0 (pq)^j = \tfrac{(pq)^j-1}{pq-1} \label{explicit expression sigma j}.
\end{align} 
Therefore, in \eqref{def D j+1} we estimate $D_j$ as follows:
\begin{align*}
D_j & = C  K^p 2^{-jp-((4+\frac{\sigma_2}{2})p+1)}   \theta_{j}^{-1} D_{j-1}^{pq} \geq C  K^p 2^{-jp-((4+\frac{\sigma_2}{2})p+1)}   (pq-1) (pq)^{-j} D_{j-1}^{pq} \\ &  = \widetilde{D} (2^p pq)^{-j}  D_{j-1}^{pq} ,
\end{align*} where $\widetilde{D}\doteq  2^{-((4+\frac{\sigma_2}{2})p+1)}  C  K^p  (pq-1) $. Repeating the same argument as in Subsection \ref{Subsection iteration argument subcrit case} (application of the logarithmic function and iterative use of the resulting inequality), we find that
\begin{align}\label{lower bound Dj}
\log D_j \geq (pq)^j \log(\widehat{D} \varepsilon),
\end{align} for $j\geq j_1\doteq \max\left\{0,\frac{\log\widetilde{D}}{\log(2^ppq)}-\frac{pq}{pq-1}\right\}$, where $\widehat{D}\doteq M (2^p pq)^{-\frac{(pq)}{(pq-1)^2}}\widetilde{D}^{\frac{1}{pq-1}}$. Combining \eqref{lower bound U j critical}, \eqref{explicit expression sigma j} and \eqref{lower bound Dj} we get
\begin{align*}
U(z) & \geq \exp\left( (pq)^j \log(\widehat{D} \varepsilon)\right)  \left(\log\left(\tfrac{z}{\ell_j R}\right)\right)^{\frac{(pq)^j-1}{pq-1} }
\end{align*} for any $j\geq j_1$ and any $z\geq \ell_{j}R$. Since $\ell_j <2$ for any $j\in\mathbb{N}$, then, from the previous inequality we get
\begin{align}
U(z) & \geq \exp\left( (pq)^j \log(\widehat{D} \varepsilon)\right)  \left(\log\left(\tfrac{z}{2R}\right)\right)^{\frac{(pq)^j-1}{pq-1} } \notag \\
 &= \exp\left( (pq)^j \left(\log(\widehat{D} \varepsilon)+\tfrac{1}{pq-1} \log\left(\log\left(\tfrac{z}{2R}\right)\right)\right) \right)  \left(\log\left(\tfrac{z}{2R}\right)\right)^{-\frac{1}{pq-1} } \notag \\
&= \exp\left( (pq)^j \log\left(\widehat{D} \varepsilon \left(\log\left(\tfrac{z}{2R}\right)\right)^{\frac{1}{pq-1}}  \right)\right)  \left(\log\left(\tfrac{z}{2R}\right)\right)^{-\frac{1}{pq-1} } \label{lower bound U final critical case}
\end{align}
for any $z\geq 2R$ and any $j\geq j_1$.

We can choose $\varepsilon_0=\varepsilon_0(n,\mu_1,\mu_2,\nu_1^2,\nu_2^2,u_0,u_1,v_0,v_1,R)>0$ sufficiently small such that $$\exp \left( (\widehat{D}\varepsilon_0)^{-(pq-1)}\right)\geq 1 .$$ Consequently, for any $\varepsilon\in (0,\varepsilon_0]$ and $z> 2R \exp \left( ( \widehat{D}\varepsilon)^{-(pq-1)}\right)$ we obtain
\begin{align*}
z\geq 2R \ \ \mbox{and} \ \ \widehat{D} \varepsilon \left(\log\left(\tfrac{z}{2R}\right)\right)^{\frac{1}{pq-1}} > 1
\end{align*} and, hence, as $j\to \infty$ in \eqref{lower bound U final critical case} we see that the lower bound for $U(z)$ is not finite. Therefore, in order to guarantee the existence of $U(z)$, it must hold the converse inequality for $z$. Finally, since we are on the characteristic $t=z+R \lesssim z$, we derive the upper bound for the  lifespan $$T(\varepsilon) \leq \exp\left(\bar{D} \varepsilon^{-(pq-1)}\right),$$ for a suitable positive constant $\bar{D}$.

Finally, in the case $ \Lambda(n+\sigma_2,q,p)=0>\Lambda(n+\sigma_1,p,q)$, by switching the role of $U$ and $V$ we get the same kind of upper bound estimate for the lifespan.

\subsubsection*{Case $\Lambda(n+\sigma_1,p,q)= \Lambda(n+\sigma_2,q,p)=0$}

In the cusp point of the critical curve we can improve the upper bound of the lifespan further. According to Remark \ref{Remark double critic case} in this case $p=\widetilde{p}(n,\sigma_1,\sigma_2)$ and $q=\widetilde{q}(n,\sigma_1,\sigma_2)$. Therefore,
\begin{align*}
-\tfrac{n-1}{2}(p-1)+\tfrac{\sigma_1}{2}-\tfrac{\sigma_2}{2}p & =  -\tfrac{n+\sigma_2-1}{2}(p-1)+\tfrac{\sigma_1}{2}-\tfrac{\sigma_2}{2} \\
& =  -\tfrac{n+\sigma_2-1}{2}\left(\tfrac{n+\sigma_1+1}{n+\sigma_2-1}-1\right)+\tfrac{\sigma_1}{2}-\tfrac{\sigma_2}{2}=-1
\end{align*} and, similarly, $-\tfrac{n-1}{2}(q-1)+\tfrac{\sigma_2}{2}-\tfrac{\sigma_1}{2}q=-1$. 

Therefore, in the case $\Lambda(n+\sigma_1,p,q)= \Lambda(n+\sigma_2,q,p)=0$ the iteration frame is simply
\begin{align}
U(z) \geq M\varepsilon +C\int_R^z  \left(R+y\right)^{-1}|V(y)|^p \, \mathrm{d}y, \label{iteration frame lower bound U double critic} \\ 
V(z) \geq N\varepsilon +K\int_R^z  \left(R+y\right)^{-1}|U(y)|^q \, \mathrm{d}y \label{iteration frame lower bound V double critic} 
\end{align} for any $z\geq R$.
Due to the special structure of \eqref{iteration frame lower bound U double critic} and \eqref{iteration frame lower bound V double critic}, in this case is not necessary to applying the slicing procedure in order to restrict step by step the domain of integration.

Thus, the first step will be to prove
\begin{align}\label{lower bound U j double critical}
U(z)\geq E_j \left(\log\left(\frac{z}{ R}\right)\right)^{\varrho_j} \qquad \mbox{for} \ \ z\geq  R,
\end{align} where as usual  $\{E_j\}_{j\in\mathbb{N}}$, $\{\varrho_j\}_{j\in\mathbb{N}}$ are suitable sequences of nonnegative real numbers. For $j=0$ we get that \eqref{lower bound U j double critical} is fulfilled provided that $E_0\doteq M\varepsilon$ and $\varrho_0\doteq 0$. We prove now the inductive step. Noticing that $R+y\leq 2y$ for $y\geq R$, if we plug \eqref{lower bound U j double critical} in \eqref{iteration frame lower bound V double critic}, then, it follows
\begin{align*}
V(z)& \geq K\int_R^z  \left(R+y\right)^{-1}|U(y)|^q \, \mathrm{d}y \geq   K E_j^q \int_R^z  \left(R+y\right)^{-1} \left(\log\left(\tfrac{y}{ R}\right)\right)^{q \varrho_j} \, \mathrm{d}y  \\
& \geq  2^{-1} K E_j^q \int_R^z  \tfrac1y \left(\log\left(\tfrac{y}{ R}\right)\right)^{q \varrho_j} \, \mathrm{d}y =  2^{-1} K E_j^q (q \varrho_j+1)^{-1} \left(\log\left(\tfrac{z}{ R}\right)\right)^{q \varrho_j+1} 
\end{align*} for any $z\geq R$. Next we use this lower bound for $V$ in \eqref{iteration frame lower bound U double critic}, obtaining
\begin{align*}
U(z) & \geq  C\int_R^z  \left(R+y\right)^{-1}|V(y)|^p \, \mathrm{d}y \\ 
& \geq 2^{-p} C K^p (q \varrho_j+1)^{-p}  E_j^{pq} \int_R^z  \left(R+y\right)^{-1}\left(\log\left(\tfrac{y}{ R}\right)\right)^{pq \varrho_j+p}  \, \mathrm{d}y \\
& \geq 2^{-(p+1)} C K^p (q \varrho_j+1)^{-p}  E_j^{pq} \int_R^z  \tfrac1y \left(\log\left(\tfrac{y}{ R}\right)\right)^{pq \varrho_j+p}  \, \mathrm{d}y \\ 
&= 2^{-(p+1)} C K^p (q \varrho_j+1)^{-p} (pq \varrho_j+p+1)^{-1}  E_j^{pq} \left(\log\left(\tfrac{z}{ R}\right)\right)^{pq \varrho_j+p+1}
\end{align*}
 for any $z\geq R$. Hence, we proved \eqref{lower bound U j double critical} for $j+1$ provided that 
 \begin{align}
 \varrho_{j+1} & \doteq p+1+  pq\,  \varrho_j , \label{def rho j+1}\\
 E_{j+1} & \doteq  2^{-(p+1)} C K^p (q \varrho_j+1)^{-p} (pq \varrho_j+p+1)^{-1}  E_j^{pq} . \label{def E j+1}
 \end{align} We determine the value of $\varrho_j$ by using \eqref{def rho j+1} iteratively
\begin{align}
\varrho_j & = p+1 + pq \varrho_{j-1}= (p+1)(1+pq)+(pq)\varrho_{j-2} = \cdots = (p+1) \sum_{k=0}^{j-1} (pq)^k +(pq)^j \varrho_0\notag \\ & = \tfrac{p+1}{pq-1}((pq)^j-1).\label{explicit expression rho j}
\end{align} Using \eqref{explicit expression rho j}, from \eqref{def E j+1} we find 
\begin{align*}
E_{j} & = 2^{-(p+1)} C K^p (q \varrho_{j-1}+1)^{-p} (pq \varrho_{j-1}+p+1)^{-1}  E_{j-1}^{pq} \\ &  \geq 2^{-(p+1)} C K^p (pq \varrho_{j-1}+p+1)^{-(p+1)}  E_{j-1}^{pq} \\ & = 2^{-(p+1)} C K^p \varrho_{j}^{-(p+1)}  E_{j-1}^{pq} \geq  \underbrace{ 2^{-(p+1)} C K^p\left(\tfrac{p+1}{pq-1}\right)^{-(p+1)}}_{\doteq \, \widetilde{E}}(pq)^{-(p+1)j}  E_{j-1}^{pq} .
\end{align*} As in the previous cases, from the inequality $E_{j}\geq  \widetilde{E} (pq)^{-(p+1)j}  E_{j-1}^{pq} $ it follows the estimate
\begin{align}\label{lower bound Ej}
E_j \geq \exp\left((pq)^j \log (\widehat{E}\varepsilon)\right)
\end{align} for any $j\geq j_2\doteq \max\{0,\frac{\log \widetilde{E}}{\log(pq)^{p+1}}-\frac{pq}{pq-1}\}$, where $\widehat{E}\doteq M (pq)^{-\frac{(pq)(p+1)}{(pq-1)^2}}\widetilde{E}^{\frac{1}{pq-1}}$. Therefore, combining \eqref{lower bound U j double critical}, \eqref{explicit expression rho j} and \eqref{lower bound Ej}, we have
\begin{align}
U(z)& \geq E_j \left(\log\left(\tfrac{z}{ R}\right)\right)^{\frac{p+1}{pq-1}((pq)^j-1)} \geq \exp\left((pq)^j \log (\widehat{E}\varepsilon)\right) \left(\log\left(\tfrac{z}{ R}\right)\right)^{\frac{p+1}{pq-1}((pq)^j-1)} \notag\\
& \geq \exp\left((pq)^j \log \left(\widehat{E}\varepsilon \left(\log\left(\tfrac{z}{ R}\right)\right)^{\frac{p+1}{pq-1}} \right)\right) \left(\log\left(\tfrac{z}{ R}\right)\right)^{-\frac{p+1}{pq-1}} \label{lower bound U final double critical case}
\end{align} for $z\geq R$ and $j\geq j_2$. We pick  $\varepsilon_0=\varepsilon_0(n,\mu_1,\mu_2,\nu_1^2,\nu_2^2,u_0,u_1,v_0,v_1,R)>0$ enough small such that $$\exp \left( (\widehat{E}\varepsilon_0)^{-\frac{pq-1}{p+1}}\right)\geq  1 .$$ Also, for any $\varepsilon\in (0,\varepsilon_0]$ and $z> R \exp \left( (\widehat{E}\varepsilon_0)^{-\frac{pq-1}{p+1}}\right)$ we find
\begin{align*}
z\geq R \ \ \mbox{and} \ \ \widehat{E}\varepsilon \left(\log\left(\tfrac{z}{ R}\right)\right)^{\frac{p+1}{pq-1}} > 1
\end{align*} and, hence, as $j\to \infty$ in \eqref{lower bound U final double critical case}  the lower bound for $U(z)$ blows up. Therefore, in order to guarantee the finiteness of $U(z)$, we have to require the opposite inequality for $z$. Moreover, we are on the characteristic $t=z+R \lesssim z$, so from the inequality $z\leq R \exp \left( (\widehat{E}\varepsilon_0)^{-\frac{pq-1}{p+1}}\right)$ we deduce  for a suitable constant $\bar{E}>0$ the upper bound for the  lifespan $$T(\varepsilon) \leq \exp\left(\bar{E} \varepsilon^{-\frac{pq-1}{p+1}}\right).$$ Finally, switching the role between $U$ and $V$ (that is, working with lower bound estimates for $V$ and applying the iteration frame \eqref{iteration frame lower bound U}-\eqref{iteration frame lower bound V} in the reverse order) we end up with the estimate
$T(\varepsilon) \leq \exp\left(E\varepsilon^{-\frac{pq-1}{q+1}}\right),$ for some constant $E>0$. Summarizing, we proved the last estimate in \eqref{upper bound lifespan WCS} too.

\section{Final remarks and open problems}

Let us compare our results with other results known in the literature. In \cite{LT18Glass} a blow-up result is proved for \eqref{semilinear CP derivative type} in the massless case by using the unbounded multiplier $m_1(t)\doteq (1+t)^\mu$ when $1<p\leq p_{\Gla}(n+2\mu)$. For the massless case (i.e. for $\nu^2=0$) in Theorem \ref{Thm blow-up derivative type nonlinearity} we proved  a blow-up result in the range $1<p\leq p_{\Gla}(n+\sigma)$, where 
\begin{align*}
\sigma =\begin{cases} 2\mu & \mbox{for} \ \ \mu\in [0,1), \\ 2 & \mbox{for}\ \ \mu\in [1,2), \\ \mu & \mbox{for}\  \ \mu\in [2,\infty). \\
\end{cases}
\end{align*} Since the Glassey exponent is decreasing with respect its argument, then, we improved the above cite result from  \cite{LT18Glass} in the case $\mu>1$ by enlarging the range of $p$ for which the nonexistence of global in time solutions can be proved under suitable assumptions for the data. Notice that for $\mu\in [0,1)$ we found exactly the same result as in  \cite{LT18Glass}.

Up to the knowledge of the author the weakly coupled system \eqref{semilinear WCS derivative type} has not been studied so far in the literature. Nonetheless, even in this case we obtained as possible candidate for the critical curve 
\begin{align*}
\max\left\{\frac{p+1}{pq-1}-\frac{n+\sigma_1-1}{2}, \frac{q+1}{pq-1}-\frac{n+\sigma_2-1}{2}\right\}=0,
\end{align*} which is a curve that presents a typical peculiarity of scale-invariant models: its branches are shifts of the branches of the critical curve for \eqref{WCS classical}. The presence of these shifts is due to the fact that we deal exactly with scale-invariant lower order terms. In the massless case, if we take into account a weaker kind of damping term, namely, scattering producing damping terms (which means that the time-dependent coefficients for the damping terms are nonnegative and summable functions), then, the critical condition for the powers in the nonlinear terms is exactly the same as in the corresponding classical not-damped case (cf. \cite{LT18Scatt,LT18Glass,LT18ComNon,WakYor18damp,PalTak19,PalTak19dt,PalTak19mix}).
Really recently, in \cite{ITW19} the blow-up dynamic for the semilinear wave equation with time-dependent and scattering producing damping and mass terms has been studied both in the subcritical and critical case.

Let us point out that we assumed throughout the paper that the parameters $\mu,\nu^2$ satisfy $\delta\geq 0$ for the single semilinear equation (respectively, $\mu_1,\mu_2,\nu_1^2,\nu_2^2$ satisfy $\delta_1,\delta_2\geq 0$ for the weakly coupled system). This assumptions are made in order to guarantee that the kernel functions defined by \eqref{def E(t,x;b,y)}, \eqref{def K0(t,x;y)} and \eqref{def K1(t,x;y)} are nonnegative functions. Moreover, in the blow-up argument we estimate from below the hypergeometric functions by positive constants. This choice is also sharp in the case $\delta>0$ when we consider an estimate from above. In the limit case $\delta=0$ though we have to include a logarithmic factor in the upper bound for the hypergeometric function in \eqref{def E(t,x;b,y)}. In any case, this is not an issue as in the global existence part one would deal in this specific case with a strict lower bound for the exponents. A similar situation is present in \cite{DabbPal18}  in the case $\delta=0$ but for study of the semilinear Cauchy problem with power nonlinearity  \eqref{semilinear CP power nonlinear}.

Clearly, in order to prove the fact that the conditions for the exponents that we get in this paper are actually sharp, the sufficient part has to be studied.

\section*{Acknowledgments}

The first author is supported by the University of Pisa, Project PRA 2018 49 and he is member of the Gruppo Nazionale per L'Analisi Matematica, la Probabilit\`{a} e le loro Applicazioni (GNAMPA) of the Instituto Nazionale di Alta Matematica (INdAM). The second author was partially supported by Zhejiang Provincial Nature Science Foundation of China under Grant No. LY18A010023.

\end{document}